\documentclass[12pt]{article}
\usepackage[english]{babel}
\usepackage{epsfig}
\usepackage{amsfonts}
\usepackage{amssymb}
\usepackage{amsmath}
\usepackage{amsthm}
\usepackage{latexsym}
\usepackage{graphicx}
\usepackage{bm}
\usepackage{tabularx}
\usepackage{booktabs} 
\usepackage{enumerate}
\usepackage{caption}
\usepackage{wrapfig}  

\usepackage{mathtools}  
\usepackage{eqnarray}
\usepackage{enumitem}

\usepackage[titletoc]{appendix}

\usepackage[numbers]{natbib}
\usepackage{comment}

\usepackage{bbm} 
\usepackage{url}
\usepackage{subfig}
\usepackage{hyperref}
\usepackage{color}
\usepackage{float}

\usepackage{hyperref}
\hypersetup{colorlinks=true, citecolor=blue, linkcolor=red, urlcolor=green}
 
\usepackage{bbm}
\usepackage{url}
\usepackage{color}
\usepackage{xcolor}
\usepackage[section]{placeins}

\theoremstyle{plain}
\newtheorem{theorem}{Theorem}[section]
\newtheorem{lemma}[theorem]{Lemma}

\newtheorem{proposition}[theorem]{Proposition}

\theoremstyle{definition}

\theoremstyle{remark}

\title{The Maximum Overlap Time in the \\ M/M/1 Queue}
\author{ 
Sergio Palomo 
\\
Systems Engineering
\\
Cornell University
\\
{sdp85@cornell.edu} \\ \\
Jamol Pender \footnote{Corresponding Author}
\\
School of Operations Research and Information Engineering
\\
Cornell University
\\
{jjp274@cornell.edu} 
}

\usepackage{natbib}
\usepackage{graphicx}
\usepackage{amsmath}

\begin{document}

\maketitle

\begin{abstract}
In this paper, we analyze the steady state maximum overlap time in the M/M/1 queue.  We derive the maximum overlap time tail distribution, its moments and the moment generating function.  We also analyze the steady state minimum overlap time of the adjacent customers and compute its moments and moment generating function.  Our results provide new insight on how customers become infected in the M/M/1 queue.   
\end{abstract}

\section{Introduction}

Overlap times have become an important metric for understanding the behavior of service systems.  Overlap times were first studied in the context of COVID-19 and infectious disease spread, see for example \citet{kang2021queueing, palomo2021measuring}.  However, since then the research on overlap times has exploded and become more important topic of interest, see for example \citet{delasay2022impacts, palomo2021overlap, hassin2021simple, hassin2023queueing, ko2022overlapping, palomo2023overlap,  ko2023number, xu2023queueing, gao2023overlap}.  In recent work by \citet{palomo2021measuring}, they show that the overlap time in the M/M/1 queue for any customer can be calculated using a simple Lindley style recursion.  Moreover, it also shows how to calculate the distribution of the overlap time of adjacent customers that are $k$ spaces apart and provides simulation results showing that the overlap time depends heavily on the service time distribution.  

In this work, we consider the M/M/1 queueing system model where the both the inter-arrival and service times are exponentially distributed. Most of the standard performance measures for this model have been studied for many years and are easily obtained thorough in the classic book \citet{shortle2018fundamentals}. However, we study the steady state maximum overlap time in this paper. 

The concept of maximum waiting time is of paramount importance, particularly in the context of studying the spread of infections within queues. It provides valuable insights into the probability of infection within the M/M/1 queue, emphasizing the individuals with whom one spends the most time or shares maximum overlap.

In order to analyze the steady state maximum overlap time, we first establish a novel connection between the steady state maximum overlap time and the steady state waiting time. Our new representation serves as a foundational tool to derive our main result.  In addition, we also provide some new results about the minimum overlap of the nearest neighbor in the queue i.e. the minimum overlap of the adjacent customers.  

As an aside of our analysis, we provide a new methodology for analyzing the steady state maximum overlap time of the single server queue.  Our methodology works as follows.  
\begin{enumerate}
    \item Find the distribution of the steady state waiting time.
    \item Find the distribution of the difference of a random service time minus a random inter-arrival time. 
    \item Use the new maximum overlap time representation in terms of the steady state waiting time to derive the distribution of the steady state maximum overlap time.  
\end{enumerate}

We provide the analysis of the maximum overlap time for the M/M/1 queue since we know the steady state waiting time distribution, our methodology also works for more complicated single server queueing models.

\subsection{Contributions of Our Work}

In this paper, we make the following contributions to the literature:
\begin{itemize}
\item We derive the steady state distribution of the maximum overlap time of any customer and show that is has an exponential tail distribution.  
\item We derive new moments and moment generating function for the maximum overlap time. 
\item We also analyze the minimum overlap time for an adjacent customer and show that its distribution is also exponential and derive its moments and moment generating function.  
\end{itemize}

\subsection{Organization of the Paper}
The remainder of the paper is organized as follows. Section \ref{Overlap_Model} introduces the overlap time and a recursion for computing.  It also provides some results about the steady state wait time of the M/M/1 queue.   In Section \ref{Max_overlap}, we present our main results for computing the tail distribution of the steady state maximum overlap time.  The moments and moment generating function are also calculated for the steady state maximum overlap time.  In Section \ref{Min_overlap}, we replicate the analytical approach employed in the previous section to investigate the minimum overlap among adjacent customers. Lastly, in Section \ref{secConc}, we offer a conclude and propose novel directions for future research.

\section{Overlaps in the Single Server Queue} \label{Overlap_Model}

In this section, we derive overlap times for any customer in the single server queue.  In particular we are interested in deriving the maximum overlap time in the single server queue.  However, because of the ordering of customers in the single server queue, the maximum overlap time for the $n^{th}$ customer can easily be computed from the waiting time of the $n^{th}$ customer or the waiting time of the next customer i.e. the $(n+1)^{th}$ customer.  Thus, in the single server queue, the maximum overlap for each customer is the maximum overlap of the adjacent customers since the service is ordered in first in first out fashion.  Suppose that we want to compute the overlap time of the $n^{th}$ and the $(n+1)^{th}$ customer, then we have the following recursion formula
\begin{eqnarray} \label{overlap_recursion}
O_{n,n+1} &=&  \max\left(D_{n} -  \sum^{n}_{i=1} A_{i}, 0  \right) .
\end{eqnarray}
where $D_n$ is defined as the departure time of the $n^{th}$ customer and $O_{n,n+j}$ is the overlap time of the $n^{th}$ and $(n+j)^{th}$ customers.  
Moreover, by manipulating the departure time of the $n^{th}$ customer we can show that the overlap between back to back customer is precisely the Lindley recursion and the waiting time of the $(n+1)^{th}$ customer i.e.
\begin{eqnarray*} 
O_{n,n+1} &=&  \max\left(D_{n} -  \sum^{n}_{i=1} A_{i}, 0  \right) \\
&=&  \max\left(W_{n} + S_n + \sum^{n-1}_{i=1} A_{i} -  \sum^{n}_{i=1} A_{i}, 0  \right) \\
&=& \max\left(W_{n} + S_n - A_{n}, 0  \right) \\
&=& W_{n+1}.
\end{eqnarray*}

Th recursion given in Equation \ref{overlap_recursion} says that the overlap time between the $n^{th}$ and the $(n+1)^{th}$ customers is exactly the maximum of zero and the departure time of the $n^{th}$ customer minus the arrival time of the $(n+1)^{th}$ customer.  Thus, if the $n^{th}$ customer departs before the $(n+1)^{th}$ customer arrives, then there is no overlap.  However, if the $n^{th}$ customer departs after the $(n+1)^{th}$ customer arrives, the overlap time is positive and is equal to the difference of the two times.  Using this idea, we can also derive the overlap time for the customer before as
\begin{eqnarray*}
O_{n-1,n} &=&  \max\left(D_{n-1} -  \sum^{n-1}_{i=1} A_{i}, 0  \right) \\
&=&  \max\left(W_{n-1} + S_{n-1} + \sum^{n-2}_{i=1} A_{i} -  \sum^{n-1}_{i=1} A_{i}, 0  \right) \\
&=& \max\left(W_{n-1} + S_{n-1} - A_{n-1}, 0  \right) \\
&=& W_{n}.
\end{eqnarray*}

In particular, we know from \citet{palomo2021measuring} that in steady state the distribution of the waiting time of the M/M/1 queue satisfies the following expression
\begin{eqnarray*}
\mathbb{P}\left( W_{\infty} > t \right) &=& \frac{\lambda}{\mu} e^{-(\mu - \lambda)t}.
\end{eqnarray*}

We will exploit these results in order to fully understand the behavior of the steady state maximum waiting time in the sequel.  

\section{The Maximum Overlap Time} \label{Max_overlap}

Here in this section, we derive our main results.  First we derive a novel representation of the steady state maximum overlap in terms of the steady state waiting time.  Then, we use the representation to derive the tail distribution of the steady state maximum overlap time.  We start with the maximum overlap time representation in terms of the waiting time.  

\begin{proposition}
Let $M_n$ be the maximum overlap time for customer $n$ in the G/G/1 queue.  Then the $M_n$ is equal to
\begin{equation}
   M_n = W_n \cdot \{ S_n < A_n \} + \left( W_n + S_n - A_{n} \right) \cdot \{ S_n \geq  A_n \}.
\end{equation}

\begin{proof}
\begin{eqnarray*} 
M_n &=& \max \left( O_{n-1,n}, O_{n,n+1} \right) \\
&=& \max \left( W_{n}, W_{n+1} \right) \\
&=& \max \left( W_{n}, \max \left( W_{n} + S_n - A_{n}, 0  \right) \right) \\
&=& \max \left( W_{n}, W_{n} + S_n - A_{n} \right) \\
&=& W_n \cdot \{ S_n < A_n \} + \left( W_n + S_n - A_{n} \right) \cdot \{ S_n \geq  A_n \}.
\end{eqnarray*} 
This completes the proof.
\end{proof}
\end{proposition}

Now that we have a representation of the maximum overlap in terms of the waiting time, we need to derive the distribution of the difference of a service time minus an inter-arrival time.  This is the second ingredient for deriving the distribution of the steady state maximum overlap time.  

\begin{lemma}
Let $X$ be an exponential random variable with rate parameter $\mu$ and $Y$ be an exponential random variable with rate parameter $\lambda$.  Then the density of $X - Y$ is equal to the following expression 
\begin{equation}
    f_{X-Y}(z) =
\begin{cases}
\frac{\lambda \mu}{\lambda + \mu} e^{\lambda z}, \quad z < 0 \\
\frac{\lambda \mu}{\lambda + \mu} e^{-\mu z}, \quad z \geq 0.
\end{cases}
\end{equation}
\begin{proof}
The proof follows from \citet{palomolearning}.
\end{proof}
\end{lemma}

Now that we have the distribution of the waiting time, the distribution of the difference between a service time and an inter-arrival time and the novel representation relating the maximum overlap to the waiting time, we show how to derive the steady state maximum overlap time distribution.  This result is obtained below.  

\begin{theorem} \label{thm_dist_max}
Suppose we have an M/M/1 queue, then the steady state maximum overlap time distribution is given by
\begin{eqnarray*} 
\mathbb{P} \left(  M_\infty > t \right) &=&  \frac{2\lambda}{\mu + \lambda}  e^{-(\mu - \lambda)t}  .
\end{eqnarray*}
\begin{proof}
Let $\mathcal{S}$ and $\mathcal{A}$ be exponential random variables with rate $\mu$ and $\lambda$ respectively.  Then we can write the maximum overlap time tail distribution as 
\begin{eqnarray*} 
\mathbb{P} \left(  M_\infty  > t \right) &=& \mathbb{P} \left(  W_\infty \cdot \{ \mathcal{S} < \mathcal{A} \} + \left( W_\infty + \mathcal{S} - \mathcal{A} \right) \cdot \{ \mathcal{S} \geq  \mathcal{A} \} > t \right) \\
&=& \mathbb{P} \left(  W_\infty > t \right) \cdot \mathbb{P} \left(  \mathcal{S} < \mathcal{A} \right) + \mathbb{P} \left(  W_\infty + \mathcal{S} - \mathcal{A} > t | \mathcal{S} \geq \mathcal{A} \right) \cdot \mathbb{P} \left(  \mathcal{S} \geq \mathcal{A} \right) \\
&=& \mathbb{P} \left(  W_\infty > t \right) \cdot \frac{\mu}{\mu + \lambda} + \mathbb{P} \left(  W_\infty + \mathcal{S} - \mathcal{A} > t | \mathcal{S} \geq \mathcal{A}  \right) \cdot \frac{\lambda}{\mu + \lambda} \\
&=& \mathbb{P} \left(  W_\infty > t \right) \cdot \frac{\mu}{\mu + \lambda} + \frac{\mathbb{P} \left(  W_\infty + \mathcal{S} - \mathcal{A} > t \cap \mathcal{S} \geq \mathcal{A}  \right)}{\mathbb{P} \left( \mathcal{S} \geq \mathcal{A}  \right)} \cdot \frac{\lambda}{\mu + \lambda} \\
&=& \frac{\lambda}{\mu + \lambda}  e^{-(\mu - \lambda)t} +  \int^{\infty}_{0} \mathbb{P} \left(  W_\infty  > t - x \right) f_{X-Y}(x) dx  \\
&=& \frac{\lambda}{\mu + \lambda}  e^{-(\mu - \lambda)t} + \int^{t}_{0} \mathbb{P} \left(  W_\infty  > t - x \right) f_{X-Y}(x) dx \\
&+& \int^{\infty}_{t} \mathbb{P} \left(  W_\infty  > t - x \right) f_{X-Y}(x) dx  \\
&=& \frac{\lambda}{\mu + \lambda}  e^{-(\mu - \lambda)t} + \frac{\lambda}{\mu} \int^{t}_{0} e^{-(\mu - \lambda)(t-x)} f_{X-Y}(x) dx \\
&+& \int^{\infty}_{t} f_{X-Y}(x) dx  \\
&=&   \frac{\lambda}{\mu + \lambda}  e^{-(\mu - \lambda)t} + \frac{\lambda}{\mu } \int^{t}_{0} e^{-(\mu - \lambda)(t-x)} \frac{\lambda \mu}{\lambda + \mu} e^{-\mu x} dx \\
&+&  \int^{\infty}_{t} \frac{\lambda \mu}{\lambda + \mu} e^{-\mu x} dx   \\
&=& \frac{2\lambda}{\mu + \lambda}  e^{-(\mu - \lambda)t}.
\end{eqnarray*}
This completes the proof.
\end{proof}
\end{theorem}

\begin{figure}[h]
\centering
\includegraphics[scale=0.5]{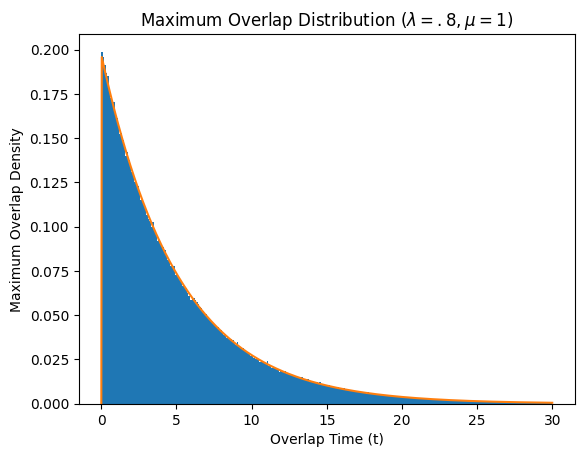}~\includegraphics[scale=0.5]{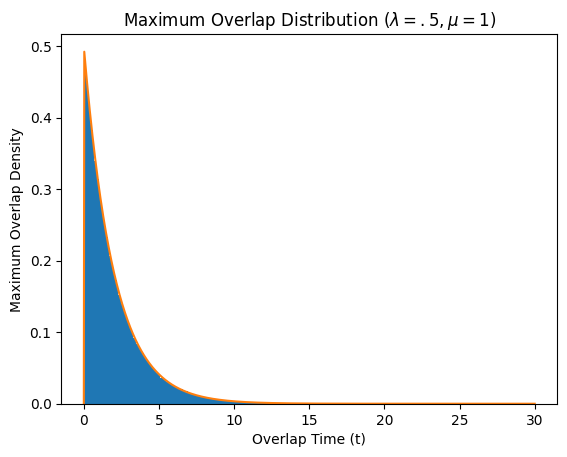}
\caption{Histogram of maximum overlap time with model parameters ($\lambda =.8, \mu = 1$) on left and ($\lambda =.5, \mu = 1$) on right.}
\label{fig:max_1}
\end{figure}

Now that we have an expression for the maximum overlap time tail distribution, we can use it to derive the moments of the maximum overlap time.  We provide this in the sequel.

\begin{theorem} \label{thm_moment}
Let $M_{\infty}$ be the steady state maximum overlap time.  Then, $p^{th}$ moment where $p\geq 1$ is equal to 
\begin{eqnarray}
    \mathbb{E} \left[ M_{\infty}^p \right] &=&  \frac{2\lambda}{\mu + \lambda} \cdot \frac{p!}{(\mu - \lambda)^p}
\end{eqnarray}

\begin{proof}
From the tail integral formula for expectation we have that
\begin{eqnarray*}
    \mathbb{E} \left[ M_{\infty}^p \right] &=&  \int^{\infty}_0 p t^{p-1} \mathbb{P} \left( M_{\infty} \geq t \right) dt \\
    &=&\int^{\infty}_0 p t^{p-1}  \frac{2\lambda}{\mu + \lambda}  e^{-(\mu - \lambda)t} dt \\
    &=& \int^{\infty}_0 p t^{p-1}  \frac{2\lambda}{\mu + \lambda}  e^{-(\mu - \lambda)t} dt \\
    &=& \frac{2\lambda p}{\mu + \lambda} \frac{1}{\mu - \lambda} \int^{\infty}_0 t^{p-1}  (\mu - \lambda)  e^{-(\mu - \lambda)t} dt \\
    &=& \frac{2\lambda p}{\mu + \lambda} \frac{1}{\mu - \lambda} \mathbb{E}_{\mu - \lambda} \left[ X^{p-1} \right] \\
    &=& \frac{2\lambda p}{\mu + \lambda} \cdot \frac{1}{\mu - \lambda} \cdot \frac{(p-1)!}{(\mu - \lambda)^{p-1}} \\
    &=&  \frac{2\lambda}{\mu + \lambda} \cdot \frac{p!}{(\mu - \lambda)^p} .
\end{eqnarray*}
This completes the proof.
\end{proof}
\end{theorem}
In particular, the moment result implies the variance of the maximum overlap time is equal to 
\begin{eqnarray}
    \mathrm{Var} \left[ M_{\infty} \right] = \frac{4 \lambda \mu}{(\lambda + \mu)^2 \cdot (\mu - \lambda)^2}.
\end{eqnarray}

\begin{theorem}
    The moment generating function of the maximum overlap time is equal to 
        \begin{eqnarray}
        \mathbb{E}\left[ e^{-\theta M_{\infty} }\right] &=& \frac{\mu - \lambda}{\mu + \lambda}  + \frac{2 \lambda}{(\mu + \lambda) \cdot (\mu +\theta - \lambda)} .
    \end{eqnarray}
    \begin{proof}
    \begin{eqnarray*}
        \mathbb{E}\left[ e^{-\theta M_{\infty} }\right] &=& \left ( 1 - \frac{2\lambda}{\mu + \lambda} \right) + \int^{\infty}_0 e^{-\theta t } \frac{2 \lambda }{\mu + \lambda} e^{-(\mu - \lambda)t} dt \\
        &=&  \frac{\mu - \lambda}{\mu + \lambda}  + \frac{2 \lambda }{ \mu + \lambda} \int^{\infty}_0  e^{-(\mu +\theta - \lambda)t} dt \\
        &=& \frac{\mu - \lambda}{\mu + \lambda}   + \frac{2 \lambda }{ (\mu + \lambda) \cdot (\mu +\theta - \lambda)} .
    \end{eqnarray*}
         This completes the proof.
    \end{proof}
    
\end{theorem}

In Figure \ref{fig:max_1}, we plot the steady state maximum overlap time with different model parameters.  It is easily observed that the distribution is exponential with same rate as the theory suggests.  

%**************************************************************************
%**************************************************************************

\section{The Minimum} \label{Min_overlap}

In addition to analyzing the maximum overlap time.  It is also interesting to analyze the minimum overlap of the adjacent customers as well.  The first step to analyzing the minimum overlap time of the adjacent customers is to derive a representation for the minimum overlap in terms of the waiting time.  This is done below.

\begin{proposition}
Let $M_n$ be the minimum overlap time for customer $n$ in the G/G/1 queue.  Then the $M_n$ is equal to
\begin{equation}
   M_n = W_n \cdot \{ S_n \geq A_n \} + \left( W_n + S_n - A_{n} \right)^+ \cdot \{ S_n <  A_n \}.
\end{equation}

\begin{proof}
\begin{eqnarray*} 
M_n^* &=& \min \left( O_{n-1,n}, O_{n,n+1} \right) \\
&=& \min \left( W_{n}, W_{n+1} \right) \\
%&=& \min \left( W_{n}, \max \left( W_{n} + S_n - A_{n}, 0  \right) \right) \\
&=& \min \left( W_{n}, (W_{n} + S_n - A_{n} )^+\right) \\
&=& W_n \cdot \{ S_n \geq A_n \} + \left( W_n + S_n - A_{n} \right)^+ \cdot \{ S_n <  A_n \}.
\end{eqnarray*} 
This completes the proof.
\end{proof}
\end{proposition}

\begin{theorem}
Suppose we have an M/M/1 queue, then the steady state minimum overlap time of the adjacent customers has the following tail distribution
\begin{eqnarray*} 
\mathbb{P} \left(  M_\infty^* > t \right) &=&   \frac{2 \lambda^2}{\mu (\mu + \lambda)}   e^{-(\mu - \lambda)t}  .
\end{eqnarray*}
\begin{proof}
Let $\mathcal{S}$ and $\mathcal{A}$ be exponential random variables with rate $\mu$ and $\lambda$ respectively.  Then we can write the minimum overlap time tail distribution as 
\begin{eqnarray*} 
\mathbb{P} \left(  M_\infty^*  > t \right) &=& \mathbb{P} \left(  W_\infty \cdot \{ \mathcal{S} \geq \mathcal{A} \} + \left( W_\infty + \mathcal{S} - \mathcal{A} \right)^+ \cdot \{ \mathcal{S} <  \mathcal{A} \} > t \right) \\
&=& \mathbb{P} \left(  W_\infty > t \right) \cdot \mathbb{P} \left(  \mathcal{S} \geq \mathcal{A} \right) + \mathbb{P} \left(  (W_\infty + \mathcal{S} - \mathcal{A})^+ > t | \mathcal{S} < \mathcal{A} \right) \cdot \mathbb{P} \left(  \mathcal{S} < \mathcal{A} \right) \\
&=& \mathbb{P} \left(  W_\infty > t \right) \cdot \frac{\lambda}{\mu + \lambda} + \mathbb{P} \left(  W_\infty + \mathcal{S} - \mathcal{A} > t | \mathcal{S} < \mathcal{A}  \right) \cdot \frac{\mu}{\mu + \lambda} \\
&=& \mathbb{P} \left(  W_\infty > t \right) \cdot \frac{\lambda}{\mu + \lambda} + \frac{\mathbb{P} \left(  W_\infty + \mathcal{S} - \mathcal{A} > t \cap \mathcal{S} < \mathcal{A}  \right)}{\mathbb{P} \left( \mathcal{S} < \mathcal{A}  \right)} \cdot \frac{\mu}{\mu + \lambda} \\
&=& \frac{\lambda^2}{\mu (\mu + \lambda)}  e^{-(\mu - \lambda)t} +  \int^{0}_{-\infty} \mathbb{P} \left(  W_\infty  > t - x \right) f_{X-Y}(x) dx  \\
&=&  \frac{\lambda^2}{\mu (\mu + \lambda)}  e^{-(\mu - \lambda)t} + \frac{\lambda}{\mu} \int^{0}_{-\infty} e^{-(\mu - \lambda)(t-x)} f_{X-Y}(x) dx \\
&=&    \frac{\lambda^2}{\mu (\mu + \lambda)}   e^{-(\mu - \lambda)t} - \frac{\lambda}{\mu } e^{-(\mu - \lambda)t} \int^{\infty}_{0} e^{-(\mu - \lambda)x} \frac{\lambda \mu}{\lambda + \mu} e^{-\lambda x} dx \\
&=&    \frac{\lambda^2}{\mu (\mu + \lambda)}   e^{-(\mu - \lambda)t} + \frac{\lambda^2 \mu}{\mu (\lambda + \mu)} e^{-(\mu - \lambda)t} \int^{\infty}_{0} e^{-\mu x}  dx \\
&=& \frac{2 \lambda^2}{\mu (\mu + \lambda)}   e^{-(\mu - \lambda)t}.
\end{eqnarray*}
This completes the proof.
\end{proof}
\end{theorem}

\begin{figure}[h]
\includegraphics[scale=0.5]{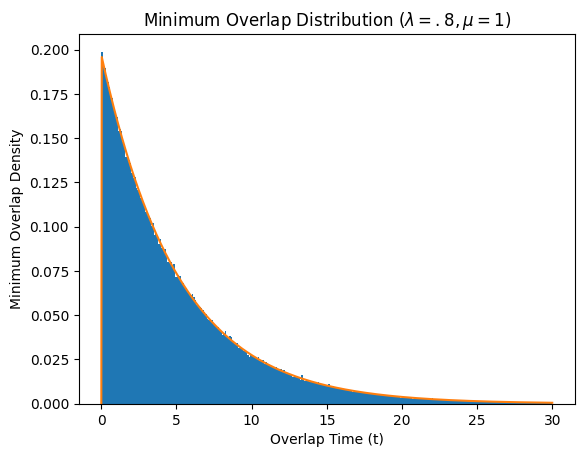}~\includegraphics[scale=0.5]{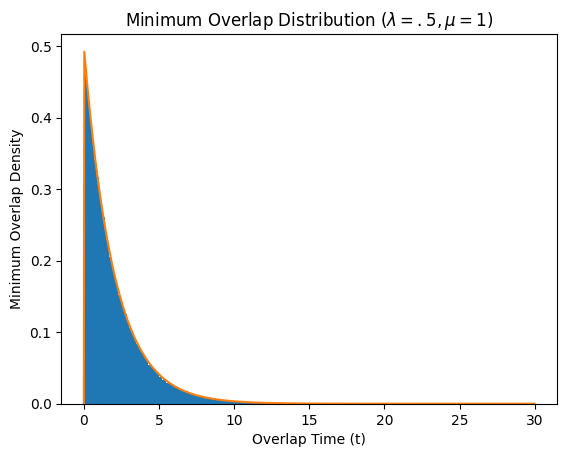}
\caption{Histogram of minimum overlap time of adjacent customers with model parameters ($\lambda =.8, \mu = 1$) on left and ($\lambda =.5, \mu = 1$) on right.}
\label{fig:min_1}
\end{figure}

\begin{theorem}
Let $M_{\infty}^*$ be the steady state minimum overlap time.  Then, the $p^{th}$ moment where $p\geq 1$ is equal to 
\begin{eqnarray}
    \mathbb{E} \left[ (M_{\infty}^*)^p \right] &=&  \frac{2\lambda^2}{ \mu( \mu + \lambda)} \cdot \frac{p!}{(\mu - \lambda)^p} .
\end{eqnarray}

\begin{proof}
The proof is identical to the proof of Theorem \ref{thm_moment} with a different constant out front.  
\end{proof}
\end{theorem}

\begin{theorem}
    The moment generating function of the minimum adjacent overlap time is equal to 
        \begin{eqnarray}
        \mathbb{E}\left[ e^{-\theta M_{\infty}^* }\right] &=& \left( 1 - \frac{2\lambda^2}{ \mu( \mu + \lambda)} \right) + \frac{2 \lambda^2}{ \mu (\mu + \lambda) \cdot (\mu +\theta - \lambda)} .
    \end{eqnarray}
    \begin{proof}
    \begin{eqnarray*}
        \mathbb{E}\left[ e^{-\theta M_{\infty}^* }\right] &=& \left( 1 - \frac{2\lambda^2}{ \mu( \mu + \lambda)} \right) + \int^{\infty}_0 e^{-\theta t } \frac{2\lambda^2}{ \mu( \mu + \lambda)} e^{-(\mu - \lambda)t} dt \\
        &=&  \left( 1 - \frac{2\lambda^2}{ \mu( \mu + \lambda)} \right)  + \frac{2\lambda^2}{ \mu( \mu + \lambda)} \int^{\infty}_0  e^{-(\mu +\theta - \lambda)t} dt \\
        &=& \left( 1 - \frac{2\lambda^2}{ \mu( \mu + \lambda)} \right) + \frac{2 \lambda^2 }{ \mu (\mu + \lambda) \cdot (\mu +\theta - \lambda)} .
    \end{eqnarray*}
         This completes the proof.
    \end{proof}
    
\end{theorem}
In Figure \ref{fig:min_1}, we plot the steady state minimum overlap time with different model parameters.  It is easily observed that the distribution is exponential with same rate as the theory suggests.

\section{Conclusion}\label{secConc}

In this paper, we analyze the maximum overlap time in the M/M/1 queue.  We find the distribution is conditionally exponential.  We also compute the moments and the moment generating function.  We also compute the distribution of the minimum overlap of the nearest customers and derive its moments and moment generating function.  These quantities are of interest as they provide new approximations for being infected in an M/M/1 queue.  The most likely customers to infect you are those that overlap with you the most.   

It would be interesting to study this maximum overlap distribution when the service times are general to understand the impact of the service time distribution.  Our proof provides a new methodology for more general distributions.  The key steps are to find the steady state waiting time distribution and the distribution of the difference of one service time and one inter-arrival time. Once these distributions are calculated, one can use our representation of the maximum overlap in terms of the waiting time to compute the maximum overlap distribution.  We hope to use this new methodology to study more complicated single server queues in the future.

%**************************************************************************
%**************************************************************************

\section*{Acknowledgements}

Jamol Pender would like to acknowledge the gracious support of the National Science Foundation DMS Award \# 2206286. 

\bibliographystyle{plainnat}
\bibliography{references}

\end{document}